# Solving a new application of asymmetric TSP by modified migrating birds optimization algorithm

Tibet Duman[a], Ekrem Duman[b,*]

[a] University of California Berkeley, Computer Science Department, Berkeley, CA, USA
{tibetduman@berkeley.edu}

[b] Ozyegin University, Industrial Engineering Department, Istanbul, Turkey
{ekrem.duman@ozyegin.edu.tr}

[*] corresponding author

*Abstract*

In this study, we first introduce a new application of the asymmetric traveling salesman problem (ATSP) which is about a small restaurant with one cook and a single stove. Once a meal has started cooking on the stove, the cook prepares the next meal on the table where the preparation time is dependent on the previous meal prepared. For the solution of this problem, besides several simple construction algorithms and a new version of the simulated annealing algorithm, we focus on enhanced versions of the recently introduced migrating birds optimization (MBO) algorithm. The original MBO algorithm might suffer from early convergence. Here we introduce several different ways of handling this problem. The extensive numerical experimentation conducted shows the superiority of the enhanced MBO over the original MBO (about 2.62 per cent) and over the simulated annealing algorithm (about 1.05 per cent).

## 1. Introduction

Recently, we came up with a new application of the TSP [1] with asymmetric and sequence dependent distance measures in a small restaurant (or, more precisely a kiosk preparing and selling different types of waffles) with one stove. After getting a number of orders, the cook prepares the materials and the raw-meal that needs to be cooked. After s/he puts it on the stove for cooking, s/he starts the preparation of the next meal. After the first meal, the preparation of the next consists of cleaning the table first and preparing the material for the next meal. The cleaning time depends on the previous meal and sometimes it can be long and sometimes short. The cooking time of each meal on the stove can be different, and before it is cooked fully, the cooking of the next meal cannot start even if it is already prepared and waiting. As detailed in the next section, this problem can be modeled as an asymmetric (ATSP).

The symmetric TSP is a very well-known problem in the literature and it is possible to see its application in many different real-life cases [1]. Likewise, ATSP is also a well-known and well-studied problem in the literature [2, 3]. The most straightforward application of ATSP is a TSP problem where distance from A to B and distance from B to A are different, perhaps due to one-way roads. Thus, while in the TSP, the roads between cities are shown by undirected edges, in the ATSP they need to be shown by directed edges [1]. Note that, while the TSP is an NP-Hard problem, the ATSP is NP-Hard in the strong sense [4]. Thus, for the solution of large sized problem instances, researchers and practitioners have always deployed heuristic algorithms.

In the broad sense, it is possible to categorize the heuristic algorithms for the TSP (and also for the ATSP) as constructive and improvement [5]. Constructive heuristics like nearest neighbor (NN) and convex-hull (CH) build solutions step by step from scratch [6]. On the other hand, improvement heuristics take a complete solution to the problem and try to improve it by small modifications. One of the small modifications that has been used in the literature is named as insertion where a randomly selected city of the given tour is removed and inserted between two other cities if it will result in a smaller cost. Another popular modification method is 2-opt where two edges of the given tour are replaced with two other edges.

Sometimes these modifications are applied in a simple manner such as trying all possible modifications and terminating when no modification results in cost reduction and sometimes, they can be applied in a more systematic way as part of metaheuristic algorithms. In metaheuristics, the new solution obtained by a modification of the current solution is usually named as a neighbor solution. It is possible to see the application of many different metaheuristic algorithms to the TSP and its variants [7-11].

We can mention the following works for SA implementations on TSP. Gent el al. [12] proposes an effective local search algorithm based on simulated annealing and greedy search techniques to solve the TSP. In order to obtain more accurate solutions, the proposed algorithm, besides following the standard simulated annealing algorithm, adopts the combination of three kinds of mutations with different probabilities during its search. Then, a greedy search technique is used to speed up the convergence rate of the proposed algorithm. Rao [13] takes up the vehicle routing problem in a supply chain network and after clustering they solve the resulting TSP instances by SA and a genetic algorithm. Rao [14] formulates the distribution problem of a FMCG company as multiple TSP and solves it by the SA algorithm. da Silva et al [15] provides a thorough study of the performance of simulated annealing in the traveling salesman problem under correlated and long tailed spatial scenarios.

Although MBO is a recently defined algorithm [16] it has been applied to many different combinatorial optimization problems. For example, Benkalai et al [17, 18], Meng et al [19], Sioud and Gagné [20], Han et al [21], Ping et al [22], Wang et al [23], Deng et al [24] and Zhang et al [25] applied it to flow shop scheduling problem while Gao and Pan [26] and Zhang et al [27] applied it to job shop scheduling problem. Xiao et al [28], Zikai [29] and Zhang et al [30] applied it to the assembly line balancing problem. Ulker and Tongur [31] solved the knapsack problem by MBO. El Aboudi and Benhlima [32] used it for feature selection in data mining. Makas and Yumuşak [33] applied it to numerical function optimization, Oz [34] to multiobjective task allocation problem, Niroomand et el [35] and Cao et al [36] to manufacturing systems, Taşpınar and Şimşir [37] to telecommunication systems, and, Tongur et al [38] to land distribution problem.

As for the MBO for the TSP, Tongur and Ülker [39] developed and compared seven different neighborhood methods for the TSP and the ATSP. They showed that the performance of MBO can be increased up to 36 per cent with the right selection of the neighborhood method. Tonyalı and Alkaya [40] applied the MBO algorithm together with two other metaheuristics (SA and ABC (artificial bee colony algorithm)) on a special variant of the TSP. They implemented and compared 10 different neighborhood methods and found out that 2-opt performed the best for MBO.

In the above-mentioned studies, researchers either directly implemented the original MBO or suggested some modifications on it (parallelization, multiple flocks, different neighborhood functions etc.) or, combined it with some other heuristics. However, none of these modifications addressed the structure of the benefit mechanism of the MBO. Benefit mechanism is the core structure of MBO which distinguishes it from the other metaheuristic methods.

Contributions of this study are threefold. First, a new application of the ATSP is introduced. Second, a small enhancement to SA and two important enhancements to MBO addressing its benefit mechanism are suggested. It turns out that all enhancements are reasonable and improve the performances of the standard versions of the algorithms. Thirdly, we implemented and compared three different neighborhood methods (2-opt, insertion and mixed) and observed that insertion performs the best on our ATSP instances.

The outline of the remainder of the paper is as follows. In section two, we give a more detailed description and the formulation of the special case we have encountered. The solution methods implemented and compared in this study are described in section 3. The enhancements we suggest for the SA and the MBO algorithms are described in section 4. The results of the numerical experimentation and the discussions about the results obtained are provided in section 5. The paper is concluded in section 6 where the paper is summarized and possible future research directions are indicated.

## 2. Problem Definition and Formulation

In our problem, at the beginning of the day the cook starts with a clean table and at the end of the day needs to leave the table clean for the next day. The cook receives a list of tele-orders from the customers in a relatively short interval of time so that there is no priority between the orders. The objective of the cook is preparing and cooking all orders and cleaning up the table in the shortest possible time.

After preparing the raw-meal of the first order, he puts it on the stove and then immediately starts the cleaning of the work table and preparation of the raw-meal of the next order. During this time, he needs to spend a negligible amount of time on the meal being cooked on the stove. In order to be able to start the cooking of the next meal, both the cooking of the current meal should be completed and the raw-meal preparation of the next meal should have been completed. Furthermore, before putting the newly prepared raw-meal on the stove, the cook cannot start the preparation of the further next raw-meals. This is because the table is too small and it can accommodate the materials for only one meal at a time. The raw-material preparation and cooking times for each meal can be different from the others but, once the orders are known these times are also known. On the other hand, the cleaning time of the table depends on what kind of order has been prepared before.

If we define the time between two orders A and B as the time that needs to pass from the starting time of the cooking of order A until the starting time of the cooking of order B, the cooking scheduling problem can be formulated as a TSP where orders correspond to cities. Time between order A and B ($t_{AB}$) can be formulated as:

$\max \{p_A, d_{AB}\}$

where;

$p_A$ = stove time of order A,

$d_{AB}$ = preparation time of order B after order A = cleaning time of the table from order A ($c_A$) + preparation time of order B ($prep_B$) tAB is not equal to tBA = max {pB, dBA} since pB and the two terms composing dBA (cB and prepA) are different and thus, our type of TSP will be an ATSP. Note that, even if stoving is not needed at all, the preparation of all orders could still be modeled as an ATSP since the table cleaning times are sequence dependent. Stoving times bring another dependency on the sequence but since both types of sequence dependencies can be calculated and put as a static cost matrix, our problem can still be regarded as a standard ATSP. We do not need to put here the mathematical model of the ATSP since it can be found anywhere but we provide some more discussion about our problem setting and the implicit and explicit assumptions we made.

We implicitly assumed that the triangle inequality is satisfied and this is a reasonable assumption in our case. It is not easy to imagine a case where $d_{AB}$ values do not satisfy the triangle inequality since making two preparations will take more time than making one preparation ($d_{AC} < d_{AB} + d_{BC}$). Then, when we think about the comparison of $t_{AC}$ with ($t_{AB} + t_{BC}$), we can argue that either $p_X$ or $d_{XY}$ will come out from $\max \{p_X, d_{XY}\}$so that if $d_{XY}$ values are in general larger than $p_X$ values, then the $d_{XY}$ values will dominate the max functions which are already satisfying the triangle inequality. Or, if $p_X$ values are usually dominating the max functions then we will be comparing a single stove time to two stove times in the comparison of $t_{AC}$ with ($t_{AB} + t_{BC}$), so, the triangle inequality will be satisfied again.

We assumed that the cook cannot start the preparation of the next order once he finished the preparation of an order and that he has to wait for the order on the stove. If this is not the case and the orders can be queued up in front of the stove while the cook prepares the next orders, the problem would get much more complicated. The problem would be a different variant of the sequence dependent TSP (SDTSP) [41] and this case is out of the scope of this study.

We also assumed that there is only one cook and one stove. In case we have multiple cooks – one stove or, one cook – multiple stoves or, multiple cooks – multiple stoves the problem formulations we would face would be different and more complicated than the one we formulated here. These cases are also left out of the scope of this current study.

## 3. Solution methods implemented

In this section we provide short descriptions of the methods that we have implemented in this study.

### *3.1.Constructive Heuristics*

Nearest Neighbor (NN):

It is one of the simplest and straightforward greedy heuristics that is used to solve the TSP. Starting from the origin each time the city that can be reached the quickest is visited. With this method, high costs can be incurred for the cities that are visited last or to come back to origin and thus, it is not a very good performing method. However, because of its simplicity, it is usually implemented to have a quick upper bound for the problem and/or to obtain an initial solution for the improvement methods.

Although NN is usually implemented for the symmetric TSP instances, here we used it for our asymmetric TSP as well where, to measure the distance to other nodes we used the $t_{AB}$ formula es explained in the previous section.

Convex-Hull:

It is a heuristic procedure that starts with a subtour consisting of the convex hull of all points to be visited. Then, at each iteration, a candidate city not on but closest to the current subtour is determined and included into the subtour by eliminating the closest edge and connecting its endpoints to the candidate point [6].

Especially for the symmetric TSPs, CH is known to be a good performing constructive method which can also be used as a starting solution for the improvement methods. Here, we implemented it assuming we have a symmetric problem (based on the $d_{AB}$ values as described above) and once the CH is obtained, we calculate its exact cost using the asymmetric cost values ($t_{AB}$).

### *3.2. Improvement Methods*

In the implementation of all improvement methods described in this section, we used the exact asymmetric cost figures.

Or-Opt:

Or-opt is usually coupled with CH (but of course can be applied to any other solution) and consists of making a systematic search of modifications that could improve the current solution. It consists of three steps [42]:

i. Starting with some first city in the given tour, consider all three consecutive cities; temporarily remove them from the tour and consider inserting them in their normal order or reverse order between any two other consecutive points in the tour (while considering the point of insertion, start with the two cities coming right after the removed three cities and proceed clockwise). Make the first insertion that yields an improvement in tour cost permanent. Continue testing other three consecutive city exchanges until the start point is reached.
ii. Repeat the procedure above for all two consecutive city exchanges.
iii. Repeat the procedure above for all single city exchanges

2-Opt

In this method, two edges of the current tour, say AB and CD, are replaced by two other edges (AC and BD) if it will bring in cost reduction. Note that, if such a change is implemented, the visiting order of the cities from B to C need to be reversed.

### *3.3. Simulated Annealing*

SA or other metaheuristic approaches try to improve a current solution or a set of current solutions by generating neighbor solutions of it and replacing it by a neighbor solution if it is better than the current solution. In SA, additionally, a worse neighbor solution can also replace the current solution probabilistically. The probability of accepting worse solutions is high at the beginning and it deteriorates as time passes. This way, it opens the way to escape from local optima. More precisely, if;

$\Delta$ = the objective value of the neighbor solution – the objective value of the current solution

Then,

If $\Delta < 0$, neighbor solution is accepted definitely

Else neighbor solution is accepted if U < $e^{-\Delta/T}$

Where U is a uniform random number between 0 and 1.

At the beginning we start with a large value for $T$, after making $L$ iterations at that temperature, we multiply it by $\alpha$ which is known as the cooling parameter and takes values close to but less than 1.00. The procedure is repeated for a predetermined number ($K$) of iterations (the number of tested neighbor solutions).

The parameters $T$, $L$ and $\alpha$ are very important for the success of SA and they can be different for each problem type and instance. Thus, the best values of them (for a given iteration limit K) should be experimentally determined before finalizing the algorithm.

### *3.4.MBO*

As opposed to SA, MBO starts with multiple ($n$) initial solutions and arranges them like on a hypothetical V shape, naming the first solution as the leader solution (bird) and the remaining as the left tail and the right tail. Each solution generates a number of its neighbors to see if there is an improvement. The number of neighbors generated and the replacement mechanism is a bit different for the leader solution and the others.

The leader solution generates $k$ neighbors and if the best of them is better than the leader solution, it is replaced by that best one. The leader solution shares $x$ best unused neighbors with its follower on the left tail and another $x$ best unused neighbors with its follower on the right tail. The other solutions generate ($k$-$x$) neighbors themselves and after combining them with the $x$ neighbors borrowed from the solution in their front, again they have $k$ neighbors to choose from. If the best of these $k$ neighbors is better, it replaces the current solution. The $x$ best unused neighbors are shared with the solution that follows.

Once neighbors are considered for all solutions till the end of the tails of the V shape, the procedure is repeated again starting from the leader solution. After $m$ repetitions, the leader changes. The leader goes to the end of one of the tails and its immediate successor in that tail becomes the new leader. The same process is repeated with the updated leader. It continues until the iteration limit $K$. More details of the MBO algorithm can be found in Duman et al [16].

As for any metaheuristic algorithm, parameter finetuning is crucial for the MBO. Here the parameters that need to be carefully studied are $n$, $k$, $x$ and $m$.

### *3.5. Neighborhood Generation*

There are many different ways of generating neighbor solutions of a given solution [39]. In our study, for both of the SA and the MBO algorithms, we implemented two types of modifications to generate neighbors: i-2-opt, ii-insertion. Insertion is the same as the third step of Or-opt where a city randomly determined is tried to be inserted into another position of the tour.

We defined three different versions of SA and MBO where the first versions implement 2-opt, second versions implement insertion and third versions implement one of these two randomly (mixed).

## 4. Enhancements on SA and MBO

### *4.1. SA*

In the description of the standard SA algorithm, it is often left unclear if the best solution found so far is kept track or not [12]. It is inherently assumed that after a large number of iterations the solution is at one of the good valleys which might include the global optimum. However, it might be the case that, although we have been on the path that would lead to the global optimum, we could have jumped to another valley since we were accepting worse solutions (see figure 1). Thus, as an enhancement to standard SA, we propose going back to the best solution found so far after 90 per cent of the iterations and spend the rest of the iterations from there on allowing downhill moves only. We name this version of SA as SA-BSF where BSF stands for best so far.

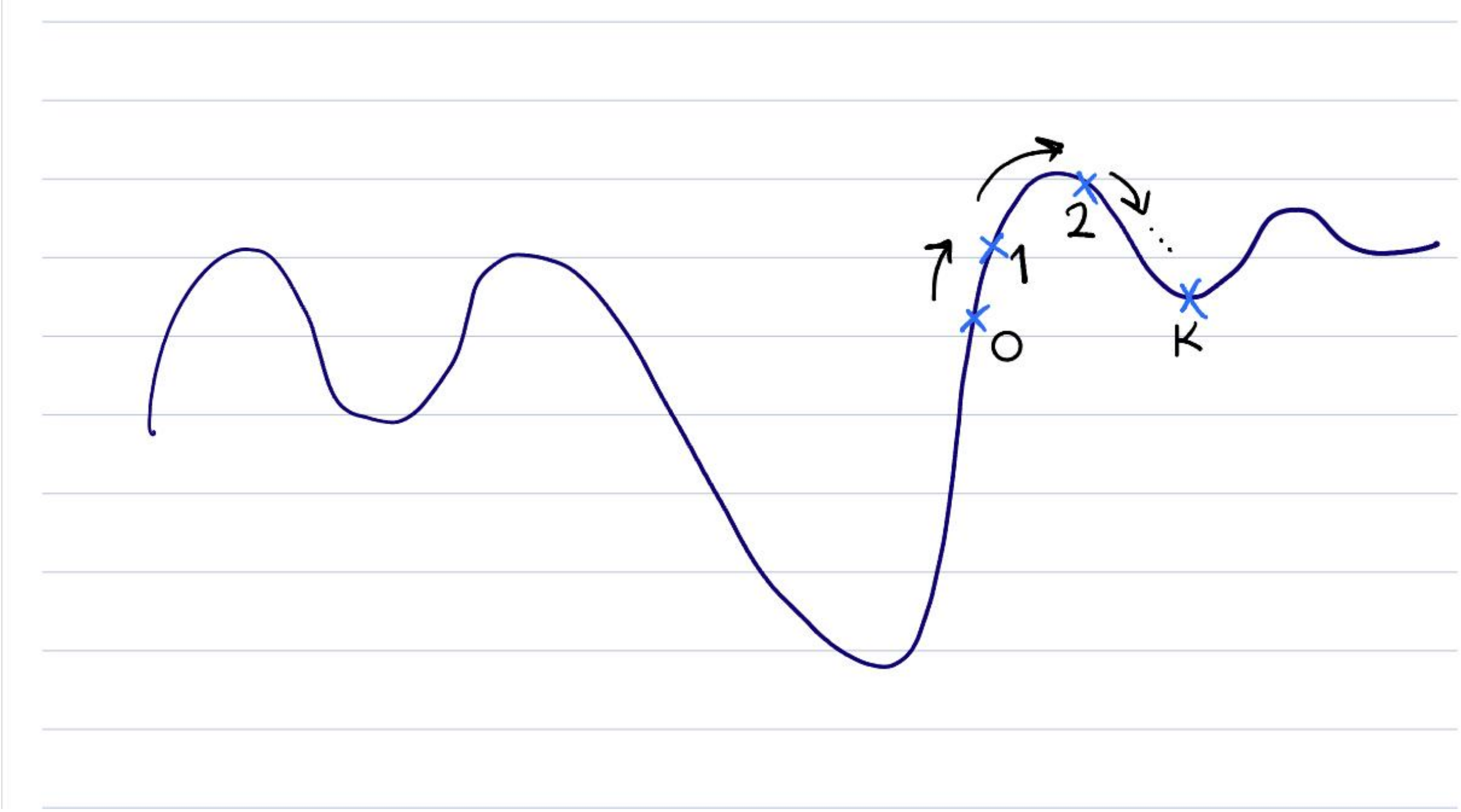

Figure 1. A possible behavior in SA for a minimization problem (0: starting solution, 1: first solution accepted, 2: second solution accepted, etc.)

Our next suggestion is about the determination of the initial temperature T in a systematic way. Since acceptance probability is dependent both on $\Delta$ and T, the value of T should be determined in accordance with $\Delta$. Obviously, the magnitude of possible $\Delta$ values change with the size of the problem instances. Thus, it is reasonable to attach the acceptance probability to some percent (say, x) worse solutions are accepted with some probability (say, y). For example, at the very first iterations of SA, it may sound reasonable to accept 10 percent worse solutions with 50 per cent probability. However, since SA in general starts with randomly generated initial solutions and the objective values of these solutions can vary a lot, accepting a modification that results in a 10 per cent worse solution of already a very bad solution could take us to very bad areas of the search space, making it difficult to bounce back to more promising areas later. Thus, we propose fixing the reference solution to a good one from which we measure 10 per cent (or, x per cent) deviation. In that regard, the solution obtained by the CH and Or-Opt couple can be a good reference.

### *4.2.MBO*

The original MBO algorithm as described above and detailed in Duman et al [16] might suffer from early convergence because of two reasons. First, a solution immediately considers (IC = immediate consideration) the best unused neighbors of the solution in the front together with its own neighbors. This may cause the solution to leave its area quickly without enough exploitation and move in the neighborhood of the solution in the front. As opposed to this, we suggest delayed consideration (DC) of the neighbors borrowed where, a solution will try to improve itself by its own neighbors and, only if it is unsuccessful, it will look at the best borrowed neighbor. This way, there is a higher probability for the exploitation of the surrounding of all solutions in the team.

Secondly, in the original MBO it is allowed to share the borrowed neighbors from the front solution with the solutions that follow (M = multi step sharing). This can be another reason for early convergence so that more solutions can quickly gather in the same neighborhood. As opposed to this, we can also consider one step sharing (S) where borrowed neighbors cannot be shared again.

In short, we propose testing the following variants of the original MBO where MBO-IC-M corresponds to the original MBO algorithm:

MBO-IC-S

MBO-IC-M

MBO-DC-S

MBO-DC-M

Pseudocode of MBO-DC-S is given below. For the convenience of the reader its parts that are different than the original MBO are underlined.

*Pseudocode of MBO-DC-S:*

1. Generate $n$ initial solutions in a random manner and place them on an hypothetical V formation arbitrarily.
2. $i$=0
3. while($i$<$K$)
   4. for ($j$=0;$j$<$m$;$j$++)
   5. Try to improve the leading solution by generating and evaluating $k$ neighbors of it.
   6. $i$=$i$+$k$
   7. for each solution $s_r$ in the flock (except leader)
   8. Try to improve $s_r$ by evaluating ($k$-$x$) neighbors of it. If not successful try to improve it by the $x$ unused best own neighbors of the solution in the front.
   9. $i$=$i$+($k$-$x$)
   10. endfor
   11. endfor
   12. Move the leader solution to the end and forward one of the solutions following it to the leader position.
13. endwhile
14. return the best solution in the flock

All algorithms implemented in this study are coded with Python 3.10, IDE: Visual Studio Code (VS Code).

## 5. Numerical Results

Since the real data of the waffle kiosk that was the source of inspiration of this study is too small and trivial to solve, we preferred to base our numerical studies on larger synthetically generated data. In this section we first explain how we have generated the problem instances to work on in this study. Then, we explain the parameter fine tuning experiments for SA and MBO which are then followed by the results obtained and discussions on them.

### *5.1. Experimentation Data*

To generate the problem instances for this study, we preferred to utilize the two-dimensional cartesian coordinate system of the classical TSP. We assumed a 20 cm by 30 cm rectangular area and we have generated either 20 (small size problems), or 50 (medium size problems), or 100 (large size problems) random points on this area. We have generated 10 instances from each of all three problem sizes. Each point is associated with a random processing (cooking) time between 2 minutes and 4 minutes determined by uniform distribution.

The travel time between two points is determined by dividing the Euclidean distance between them with the speed parameter. To determine the speed value, we conducted some experiments so that approximately half of the time distances between two points is determined by the cooking time. For this we needed a solution for the TSP and we assumed we could use CH since it is fast and can produce fairly good solutions. Obviously, such kinds of balanced problems are more challenging.

The random data generation Phyton code and the data itself can be found at the following url:
https://github.com/tibetduman/cook-shceduling-problem/tree/main/problems

### *5.2. Results and Discussion*

Before giving out the results, we would like to give some information on the parameter fine tuning experiments we have made. For SA (and also for MBO), we preferred to set the iteration limit K (for the number of neighbor solutions to be generated) to five times the cube of the problem size. Thus, K was taken as equal to 40.000, 625.000 and 5.000.000 for the small, medium and large size problems respectively. To determine the initial temperature, we tested three alternatives: accepting 10 per cent worse solutions with 75, 50 and 25 per cent probability. The final temperature is determined as the one corresponding to accepting 1 per cent worse solutions with 0.0001 probability. To determine the number of neighbors generated at each temperature we wanted to fix the number of cooling steps to be applied across all three problem sizes. In this regard, the neighbor numbers tested were 10 and 20 for small problems, 156 and 312 for medium problems and 1250 and 2500 for large problems. After we noticed that large *L* values are resulting in escaping from good valleys which are difficult to bounce back later on, for large problems we tested 312 as well. As a result of the extensive experimentation, the best initial starting temperature corresponded to accepting 10 per cent worse solutions for all problem sizes. The best *L* values were 10, 312 and 312 and the corresponding best alfa values were 0.995, 0.986 and 0.9987 for the small, medium and large problems respectively.

As for the MBO algorithm, at the beginning we wanted to fix the value of *x* to 1 for two reasons. First, in the original study of Duman et al (2012) $x = 1$ was supported by the experiments. Second, for single step sharing algorithms (MBO-DC-S and MBO-IC-S) *x* greater than 1 is meaningless. For the parameter *m* we tested different values of 1, 2, 5, 10 and 20 and since the algorithm was insensitive to these values, we preferred to keep it as equal to 1 for simplicity. Then, experimentation for small problems included *n* (number of birds) values of 11, 21, 31, 51 and *k* (number of neighbors) of 5, 7, 11. The conclusion we reached from our experimentation was that when the bird number is 11, regardless of *K*, in 40.000 neighbor generations all of the birds very easily converged to the same solution, and often enough, this solution was not a good performing one. This stemmed from the fact that with 11 birds, we do not have enough samples from the solution space to explore sufficiently. When the bird number was 31 and higher like 51, we noticed that the birds did not converge and needed a very high corresponding *k* value so that each solution is exploited sufficiently. However, we found that when the number of birds was 21 it seemed to give us enough random starting positions from the solution space and just enough replacements so that each bird exploited its neighborhood and improved itself as much as it could. When it came to *k*, the number of neighbors each bird generated, we noticed that when *k* was 5, 40.000 neighbors usually was not enough for the birds to reach their corresponding local optima, we found 7 to be sufficient for them to reach their local optima in the permitted neighbor count. When testing higher values of *k*, like 11, we noticed that it often gave almost no improvement except for the setups with very high numbers of birds, which had difficulty to converge. So, we decided that $n = 21$ and $k = 7$ as the most reasonable parameters for this problem size. Similar experimentation with similar arguments made for medium and large problems ended up with $n = 101$, $k = 7$ for medium problems and $n = 501$, $k = 11$ for large problems where in both cases $m = x = 1$ was preferred again.

The results (total operation times in seconds) of the small, medium and large problems are given in Tables 1, 2 and 3 respectively. Each of the small problems are run 10 times and the best, the worst and the average of them are tabulated. Each of the medium problems are run three times and the best and average of them are tabulated. On the other hand, the large problems are run only one time since the run times are quite large (in Table 3, MBO is shortened to M to save space). In each row, the solution with the lowest objective value is bolded. All problems are solved optimally by Gurobi and the result is shown in the *Solver* column. The average performances of the algorithms in terms of per cent deviation from the optimum solution are calculated and shown in Table 4.

For the small problems, in terms of the averages, the best results are obtained by the MBO-DC-S algorithm with insertion neighborhood, followed by the MBO-DC-M with insertion neighborhood, MBO-DC-S algorithm with mixed neighborhood, MBO-DC-M algorithm with mixed neighborhood and SA-BSF with mixed and insertion neighborhoods. In terms of the best of the 10 runs, the best algorithms having the same score are MBO-DC-S and MBO-DC-M algorithms with insertion neighborhood, MBO-DC-M and MBO-IC-S algorithms with mixed neighborhood. These are followed by the MBO-IC-M algorithm with insertion neighborhood and SA-BSF algorithm with insertion neighborhood that are having a slightly worse score. MBO-DC-S with insertion and mixed neighborhood and MBO-DC-M with insertion neighborhood were able to find the optimum solution for all 10 ten problems.

For the medium problems, in terms of the averages, the best results are obtained by the MBO-DC-S algorithm with insertion neighborhood, followed by the SA-BSF algorithm with insertion neighborhood, and later by the MBO-IC-S algorithm with insertion neighborhood, SA algorithm with insertion neighborhood, MBO-DC-M and MBO-DC-S algorithms with mixed neighborhood. The best solutions are obtained by the MBO-DC-S algorithm with insertion neighborhood, followed by the SA-BSF algorithm with insertion and mixed neighborhoods. MBO-DC-S algorithm with insertion neighborhood was able to obtain the optimum solution of 3 out of 10 problems.

For the large algorithms, it is not easy to arrive at a significant conclusion since the problems were run only once. However, from this it can be seen that SA and SA-BSF algorithms with insertion neighbor generation are performing well followed by SA algorithms with mixed and MBO algorithms with mixed and insertion neighborhoods. From other analyses not tabulated here we know that the iteration limit of 5 million neighbors was not enough for convergence especially for the MBO algorithms.

Table 1a. Results of Simulated Annealing for small problems (K=40.000).

| | | | | | 2-opt | | | | | | Insertion | | | | | | Mixed | | | | | |
|---|---|---|---|---|---|---|---|---|---|---|---|---|---|---|---|---|---|---|---|---|---|---|
| | | | | | SA | | | SA-BSF | | | SA | | | SA-BSF | | | SA | | | SA-BSF | | |
| Problem | NN | CH | CHOrOpt | Solver | Best | Avg | Worst | Best | Avg | Worst | Best | Avg | Worst | Best | Avg | Worst | Best | Avg | Worst | Best | Avg | Worst |
| S1 | 85.54 | 92.53 | 91.61 | 83.45 | **83.45** | 88.62 | 102.2 | **83.45** | 89.04 | 99.15 | **83.45** | 83.93 | 84.82 | **83.45** | 83.82 | 84.69 | **83.45** | 84.47 | 87.1 | **83.45** | 83.69 | 83.97 |
| S2 | 98.61 | 83.93 | 83.93 | 79.34 | 80.08 | 85.48 | 94.1 | 80.08 | 87 | 105.7 | **79.34** | 81.97 | 85.35 | **79.34** | 81.07 | 86.84 | **79.34** | 81.42 | 84.15 | **79.34** | 80.91 | 85.8 |
| S3 | 92.6 | 103.6 | 98.81 | 81.26 | **81.26** | 89.9 | 95.73 | **81.26** | 88.45 | 96.18 | **81.26** | 82.69 | 85.4 | **81.26** | 82.99 | 89.01 | **81.26** | 83.46 | 88.88 | **81.26** | 84.19 | 88.54 |
| S4 | 101 | 93.86 | 92.69 | 83.36 | 84.29 | 87.35 | 92.41 | 83.8 | 87.17 | 91.44 | 83.71 | 85.39 | 87.62 | 83.71 | 84.65 | 86.99 | 83.71 | 85.04 | 87.8 | **83.36** | 84.37 | 87 |
| S5 | 82.81 | 92.76 | 89.9 | 77.47 | 78.59 | 82.55 | 94.46 | **77.47** | 80.76 | 84.26 | **77.47** | 78.01 | 79.6 | **77.47** | 78.62 | 79.82 | **77.47** | 79.01 | 80.9 | **77.47** | 78.31 | 83.69 |
| S6 | 98.24 | 87.8 | 87.52 | 82 | 82.37 | 86.25 | 103.8 | **81.99** | 86.04 | 94.59 | **81.99** | 82.82 | 83.96 | **81.99** | 82.9 | 83.96 | **81.99** | 83.72 | 90.24 | **81.99** | 83.5 | 91.16 |
| S7 | 86.48 | 96.15 | 95.35 | 80.65 | **80.65** | 85.1 | 93.23 | 81.68 | 85.77 | 92.7 | **80.65** | 82.33 | 83.71 | **80.65** | 81.83 | 83.08 | **80.65** | 82.01 | 82.82 | **80.65** | 81.72 | 82.82 |
| S8 | 87.21 | 90 | 89.9 | 75.19 | 76.45 | 83.84 | 89.13 | **75.19** | 79.03 | 82.73 | **75.19** | 77.7 | 83.68 | **75.19** | 77.02 | 83.4 | **75.19** | 77.23 | 82.12 | **75.19** | 76.01 | 77.12 |
| S9 | 100 | 87.67 | 86.82 | 82.55 | **82.55** | 86.25 | 93.23 | 84.22 | 85.74 | 90.75 | **82.55** | 84.18 | 85.94 | 82.95 | 84.93 | 88.42 | 82.95 | 84.32 | 85.95 | 82.95 | 84.07 | 84.91 |
| S10 | 94.6 | 89.79 | 89.77 | 86.05 | 87.93 | 90.15 | 102.1 | 86.51 | 89.38 | 98.66 | **86.05** | 87.52 | 89.25 | **86.05** | 86.77 | 87.97 | **86.05** | 87.19 | 88.41 | **86.05** | 87.7 | 89.45 |
| Average | 92.71 | 91.81 | 90.63 | 81.13 | 81.76 | **86.55** | 96.05 | 81.57 | **85.84** | 93.62 | 81.17 | **82.65** | 84.93 | 81.21 | **82.46** | 85.42 | 81.21 | **82.79** | 85.84 | 81.17 | **82.45** | 85.45 |

Table 1b. Results of MBO for small problems (K=40.000).

| | 2-opt | | | | | | | | | | | | Insertion | | | | | | | | | | | |
|---|---|---|---|---|---|---|---|---|---|---|---|---|---|---|---|---|---|---|---|---|---|---|---|---|
| | MBO-IC-M | | | MBO-IC-S | | | MBO-DC-M | | | MBO-DC-S | | | MBO-IC-M | | | MBO-IC-S | | | MBO-DC-M | | | MBO-DC-S | | |
| Problem | Best | Avg | Worst | Best | Avg | Worst | Best | Avg | Worst | Best | Avg | Worst | Best | Avg | Worst | Best | Avg | Worst | Best | Avg | Worst | Best | Avg | Worst |
| S1 | 83.76 | 90.43 | 104.1 | **83.45** | 86.82 | 92.43 | **83.45** | 84.57 | 86.38 | **83.45** | 84.42 | 86.67 | **83.45** | 84.87 | 91.17 | **83.45** | 84.22 | 85.78 | **83.45** | 83.68 | 83.97 | **83.45** | 83.75 | 83.97 |
| S2 | **79.34** | 84.33 | 97.47 | 80.01 | 86.11 | 94.66 | **79.34** | 82.36 | 86.75 | **79.34** | 81.89 | 85.31 | 82.1 | 83.06 | 85.83 | **79.34** | 81.68 | 88.98 | **79.34** | 79.7 | 82.19 | **79.34** | **79.34** | **79.34** |
| S3 | **81.26** | 90.97 | 97.13 | 85.65 | 88.76 | 93.43 | 82.4 | 86.06 | 89.08 | **81.26** | 84.18 | 88.47 | **81.26** | 82.46 | 85.59 | **81.26** | 84.32 | 88.69 | **81.26** | 82.17 | 85.18 | **81.26** | 81.52 | 82.59 |
| S4 | 84.69 | 87.6 | 91.76 | 85.61 | 87.65 | 89.38 | 83.8 | 85.45 | 86.92 | 83.8 | 85.01 | 86.35 | **83.36** | 84.85 | 87 | **83.36** | 84.2 | 85.96 | **83.36** | 83.61 | 84.29 | **83.36** | 83.92 | 84.71 |
| S5 | 78.59 | 82.84 | 88.35 | 77.5 | 81.08 | 86.39 | **77.47** | 79.27 | 81.49 | 77.91 | 78.98 | 80.07 | **77.47** | 78.69 | 81.35 | **77.47** | 78.89 | 80.83 | **77.47** | 77.68 | 78.59 | **77.47** | 77.63 | 78.59 |
| S6 | **81.99** | 86.79 | 95.51 | **81.99** | 86.1 | 92.38 | **81.99** | 83.46 | 87.47 | 91.99 | 82.93 | 83.77 | **81.99** | 83.64 | 91.82 | **81.99** | 83.69 | 92.23 | **81.99** | 82.39 | 82.85 | **81.99** | 82.3 | 82.72 |
| S7 | 82.11 | 85.5 | 90.76 | 81.69 | 84.31 | 90.25 | **80.65** | 82.13 | 83.15 | 80.87 | 82.26 | 83.44 | **80.65** | 81.73 | 83.25 | **80.65** | 82.17 | 83.64 | **80.65** | 80.9 | 82.05 | **80.65** | 81.09 | 82.05 |
| S8 | **75.19** | 81.28 | 87.11 | **75.19** | 81.1 | 89.99 | **75.19** | 78.96 | 84.47 | **75.19** | 77.23 | 82.4 | **75.19** | 77.85 | 87.35 | **75.19** | 76.03 | 82.33 | **75.19** | 75.44 | 76.45 | **75.19** | 75.37 | 76.45 |
| S9 | 84.71 | 86.23 | 88.95 | 83.77 | 85.25 | 88.51 | 83.3 | 85.2 | 86.19 | 82.98 | 84.27 | 85.17 | 82.84 | 83.89 | 85.88 | 82.95 | 84.28 | 87.53 | **82.55** | 83.18 | 84.35 | **82.55** | 83.25 | 84.08 |
| S10 | **86.05** | 88.76 | 93.94 | 86.51 | 90.29 | 95.37 | 86.07 | 87.84 | 91.15 | **86.05** | 86.64 | 87.13 | **86.05** | 88.17 | 95.35 | **86.05** | 86.98 | 88.63 | **86.05** | 86.31 | 86.76 | **86.05** | 86.27 | 86.76 |
| Average | 81.77 | **86.47** | 93.51 | 82.14 | **85.75** | 91.28 | 81.37 | **83.53** | 86.31 | 82.28 | **82.78** | 84.88 | 81.44 | **82.92** | 87.46 | 81.17 | **82.65** | 86.46 | 81.13 | **81.51** | 82.67 | 81.13 | **81.44** | 82.13 |

Table 1c. Results of MBO for small problems continued (K=40.000).

| | Mixed | | | | | | | | | | | |
|---|---|---|---|---|---|---|---|---|---|---|---|---|
| | MBO-IC-M | | | MBO-IC-S | | | MBO-DC-M | | | MBO-DC-S | | |
| Problem | Best | Avg | Worst | Best | Avg | Worst | Best | Avg | Worst | Best | Avg | Worst |
| S1 | **83.45** | 84.51 | 87.82 | **83.45** | 84.27 | 87.46 | **83.45** | 83.77 | 83.97 | **83.45** | 83.73 | 83.97 |
| S2 | **79.34** | 80.52 | 82.45 | **79.34** | 81.84 | 84.58 | **79.34** | 80.09 | 82.99 | **79.34** | 80.14 | 82.45 |
| S3 | **81.26** | 83.44 | 88.12 | **81.26** | 83.79 | 88.28 | **81.26** | 82.31 | 85.4 | **81.26** | 82.44 | 84.42 |
| S4 | **83.36** | 84.4 | 86.78 | **83.36** | 85.07 | 87 | **83.36** | 83.65 | 84.85 | **83.36** | 83.94 | 85.87 |
| S5 | 77.98 | 79.29 | 80.26 | **77.47** | 79.42 | 86.99 | **77.47** | 78.26 | 79.6 | **77.47** | 78.11 | 80.24 |
| S6 | **81.99** | 83.01 | 84.62 | **81.99** | 83.43 | 90.44 | **81.99** | 82.26 | 82.51 | **81.99** | 82.17 | 82.37 |
| S7 | **80.65** | 82.61 | 85.96 | **80.65** | 82.11 | 89.78 | **80.65** | 81.11 | 82.08 | **80.65** | 81.28 | 82.25 |
| S8 | **75.19** | 78.16 | 87.31 | **75.19** | 76.26 | 82.09 | **75.19** | 75.32 | 84.76 | **75.19** | 75.35 | 76.04 |
| S9 | **82.55** | 84.31 | 85.16 | **82.55** | 83.82 | 85.4 | 82.84 | 83.69 | 84.76 | **82.55** | 83.29 | 84.32 |
| S10 | **86.05** | 87.77 | 98.16 | **86.05** | 86.55 | 88.41 | **86.05** | 86.49 | 88.86 | **86.05** | 86.23 | 86.76 |
| Average | 81.18 | **82.80** | 86.66 | 81.13 | **82.66** | 87.04 | 81.16 | **81.70** | 83.98 | 81.13 | **81.67** | 82.87 |

Table 2a. Results (costs) of Simulated Annealing for medium problems (K=625.000).

| | | | | | 2-opt | | | | Insertion | | | | Mixed | | | |
|---|---|---|---|---|---|---|---|---|---|---|---|---|---|---|---|---|
| | | | | | SA | | SA-BSF | | SA | | SA-BSF | | SA | | SA-BSF | |
| Problem | NN | CH | CHOrOpt | Solver | Best | Avg | Best | Avg | Best | Avg | Best | Avg | Best | Avg | Best | Avg |
| M1 | 248.9 | 205 | 204.28 | 192.3 | 204 | 214.8 | 204.4 | 209.5 | 194.8 | 196.5 | 195.1 | 195.9 | 193.6 | 196.4 | 193 | 196.3 |
| M2 | 246.9 | 225.2 | 222.61 | 194.8 | 199.6 | 212.9 | 202.5 | 210 | **197** | 201.3 | 200.3 | 202 | 197.4 | 199.1 | 198.8 | 200.6 |
| M3 | 234.5 | 221.7 | 216.71 | 189.3 | 197.4 | 204.2 | 202.5 | 207.4 | 192.5 | 194.6 | **190.2** | 193.6 | 191.3 | 196.1 | 192.5 | 196.2 |
| M4 | 241 | 199.2 | 197.9 | 190.6 | 193.4 | 210 | 197.2 | 211.5 | 192.9 | 199.4 | 194.7 | 198.1 | 195 | 202 | 194.3 | 199.4 |
| M5 | 237.1 | 218.6 | 217.94 | 194 | 214.6 | 220.8 | 204.7 | 217 | 195.5 | 200.2 | 194.6 | 198.7 | 194.6 | 206.8 | 194.6 | 199.9 |
| M6 | 213.3 | 220.2 | 216.14 | 187.2 | 193.3 | 207.8 | 201.1 | 212.4 | 187.3 | 189.1 | 187.4 | 191 | **187.2** | 193.1 | 188.6 | 192.6 |
| M7 | 259.1 | 232.8 | 230.52 | 206.3 | 214.2 | 221 | 211.2 | 224.4 | 209.2 | 213 | 209 | 212.4 | 208.8 | 212.5 | 208.9 | 211.9 |
| M8 | 220.6 | 201.4 | 198.37 | 191.3 | 196.6 | 203.1 | 198.4 | 205.5 | 193.3 | 194.6 | **191.3** | 193 | 192.8 | 195.1 | 193.1 | 194.4 |
| M9 | 235 | 221.4 | 214.02 | 189.7 | 207.2 | 219.1 | 199.8 | 211.8 | 191.6 | 198.4 | 192.1 | 195.8 | 192 | 194.5 | 191.9 | 195.5 |
| M10 | 224.4 | 232 | 228.97 | 186.9 | 200.4 | 208.9 | 199.5 | 213.6 | 192.2 | 196.3 | **186.9** | 192.6 | 189.4 | 193.1 | 190.7 | 195.3 |
| Average | 236.06 | 217.75 | 214.75 | 192.24 | 202.06 | **212.25** | 202.13 | **212.31** | 194.62 | **198.33** | 194.16 | **197.29** | 194.22 | **198.86** | 194.63 | **198.19** |

Table 2b. Results (costs) of MBO for medium problems (K=625.000).

| | 2-opt | | | | | | | | Insertion | | | | | | | | Mixed | | | | | | | |
|---|---|---|---|---|---|---|---|---|---|---|---|---|---|---|---|---|---|---|---|---|---|---|---|---|
| | MBO-IC-M | | MBO-IC-S | | MBO-DC-M | | MBO-DC-S | | MBO-IC-M | | MBO-IC-S | | MBO-DC-M | | MBO-DC-S | | MBO-IC-M | | MBO-IC-S | | MBO-DC-M | | MBO-DC-S | |
| Problem | Best | Avg | Best | Avg | Best | Avg | Best | Avg | Best | Avg | Best | Avg | Best | Avg | Best | Avg | Best | Avg | Best | Avg | Best | Avg | Best | Avg |
| M1 | 214.5 | 218.9 | 212.7 | 217.2 | 212.6 | 220 | 213.1 | 216.3 | 196.6 | 203.5 | **192.3** | 193.5 | 194.4 | 198.2 | 194 | 195.2 | 196.9 | 200.1 | 198.4 | 203.9 | 193.6 | 196.6 | 195.1 | 198.1 |
| M2 | 218 | 227.4 | 213.1 | 216.3 | 211.6 | 215.2 | 210.5 | 215.4 | 199 | 202.5 | 198.2 | 200.7 | 199.6 | 202 | 198.2 | 199.7 | 201.9 | 203 | 203.7 | 205.5 | 200 | 203.8 | 197.1 | 201.4 |
| M3 | 199.2 | 205.8 | 205.7 | 208.7 | 201.9 | 203.3 | 208.9 | 211.8 | 195.8 | 198.3 | 191 | 193.4 | 196.2 | 198.3 | 193.9 | 195.2 | 196.7 | 198.7 | 193.2 | 195.6 | 192.1 | 194.4 | 195 | 195.6 |
| M4 | 219.3 | 225.2 | 217 | 220.7 | 222.1 | 224.2 | 207.8 | 214.5 | 201.5 | 203 | 194.1 | 199.7 | 195.4 | 198.4 | **191.6** | 193 | 198.5 | 201.7 | 193.9 | 197.1 | 195.6 | 196.8 | 199.6 | 201.7 |
| M5 | 221.8 | 234.8 | 219.2 | 226 | 225 | 228.8 | 210.7 | 218.6 | 200 | 206.7 | 203.9 | 209.2 | 195.5 | 204.7 | **194** | 198.4 | 210.3 | 218.4 | 194 | 211.8 | 194.6 | 204.6 | 195.7 | 203.8 |
| M6 | 221.3 | 229.4 | 219.5 | 222.1 | 203.7 | 210.7 | 200.7 | 205.2 | 191.2 | 195.1 | 190.5 | 194.6 | 187.5 | 190.1 | 187.5 | 189 | 188.6 | 194.9 | 188.6 | 191.2 | 187.7 | 190.4 | 190.6 | 192 |
| M7 | 223 | 230.4 | 227 | 236.2 | 224 | 227.5 | 223.3 | 224.8 | 214.2 | 215 | 210.3 | 212 | 212.7 | 217.5 | **206.3** | 208.7 | 208.8 | 211.9 | 209.1 | 212.2 | 210.6 | 212.7 | 206.8 | 208.5 |
| M8 | 206.2 | 217.5 | 206.7 | 209.9 | 199 | 205.7 | 197 | 204.2 | 194.2 | 200.2 | 193.4 | 194.2 | 194.2 | 196.7 | 192.2 | 194.2 | 194 | 197.1 | 195.4 | 196.6 | 194.2 | 194.5 | 192.2 | 194.1 |
| M9 | 215.5 | 236.6 | 209.8 | 215.3 | 211.1 | 214.2 | 211.5 | 219.5 | 197.5 | 199.4 | **189.9** | 193.5 | 192.3 | 195.6 | 191.6 | 194.7 | 198.9 | 202.9 | 197.7 | 199 | 198.3 | 199.9 | 191.6 | 193.7 |
| M10 | 201.5 | 210.9 | 209 | 214.3 | 198.7 | 207.9 | 204.8 | 210 | 187.7 | 200.2 | 190.3 | 192.1 | 188.4 | 193.7 | 187.7 | 190.4 | 192.3 | 198.7 | 189.2 | 196 | 189.3 | 191.2 | 193.4 | 195.9 |
| Average | 214.04 | **223.69** | 213.96 | **218.67** | 210.95 | **215.76** | 208.82 | **214.02** | 197.78 | **202.38** | 195.39 | **198.29** | 195.61 | **199.50** | 193.70 | **195.84** | 198.68 | **202.73** | 196.32 | **200.90** | 195.60 | **198.50** | 195.71 | **198.48** |

Table 3. Results (costs) of Simulated Annealing for large problems (K=5.000.000).

| | | | | | 2-opt | | Insertion | | Mixed | | 2-opt | | | | Insertion | | | | Mixed | | | |
|---|---|---|---|---|---|---|---|---|---|---|---|---|---|---|---|---|---|---|---|---|---|---|
| Problem | NN | CH | CHOrOpt | Solver | SA | SA-BSF | SA | SA-BSF | SA | SA-BSF | M-IC-M | M-IC-S | M-DC-M | M-DC-S | M-IC-M | M-IC-S | M-DC-M | M-DC-S | M-IC-M | M-IC-S | M-DC-M | M-DC-S |
| L1 | 415.91 | 412.91 | 401.74 | 355.21 | 385.98 | 404.56 | 365.5 | 372.9 | 365.15 | 379.56 | 423.6 | 411.84 | 409.03 | 405.01 | 375.13 | 366.06 | **364.37** | 371.58 | 371.07 | 375.53 | 371 | 366.64 |
| L2 | 405.65 | 388.12 | 386.35 | 350.19 | 398.23 | 390.96 | **358.15** | 362.55 | 364.9 | 386.7 | 436.16 | 439.07 | 441.25 | 420.76 | 370.75 | 363.83 | 395.01 | 371.6 | 378.99 | 366.53 | 380.02 | 381.05 |
| L3 | 450.86 | 396.62 | 390.06 | 353.07 | 392.49 | 411.83 | **360.57** | 365.38 | 378.37 | 405.89 | 445.22 | 424.67 | 441 | 418.23 | 381.91 | 376.36 | 378.21 | 375.34 | 364.94 | 374.24 | 383.49 | 387.5 |
| L4 | 407.66 | 386.62 | 384.33 | 343.63 | 385.85 | 367.27 | **350.31** | 368.28 | 359.25 | 405.59 | 441.41 | 418.39 | 437.57 | 413.65 | 374.21 | 367.31 | 373.77 | 366.98 | 363.43 | 370.4 | 371.61 | 368.96 |
| L5 | 420.5 | 392.72 | 386.07 | 343.57 | 384.99 | 372.12 | 358.49 | **356.71** | 359.19 | 365.31 | 425.8 | 425.5 | 446.45 | 402.83 | 372.39 | 361.73 | 366.11 | 365.87 | 371.83 | 368.07 | 376.69 | 370.56 |
| L6 | 420.13 | 404.87 | 400.59 | 360.79 | 391.47 | 395.52 | 375.99 | 382 | 378.96 | 415.61 | 432.8 | 425.73 | 420.16 | 433.6 | **375.36** | 384.06 | 379.01 | 384.24 | 391.96 | 382.14 | 382.68 | 378.04 |
| L7 | 438.17 | 397.03 | 390.23 | 359.53 | 391.04 | 397.29 | 378.49 | 374.75 | 380.5 | 385.24 | 423.99 | 419.71 | 436.42 | 419.04 | 386.65 | 384.35 | 375.9 | **371.59** | 390.97 | 385.38 | 379.72 | 384.79 |
| L8 | 465.15 | 386.36 | 381.7 | 353.33 | 394.02 | 390.34 | **364.29** | 366.85 | 368.63 | 396.41 | 451.62 | 410.25 | 428.42 | 412.67 | 372.03 | 379.2 | 378.67 | 367.46 | 389.71 | 380.42 | 392.98 | 375.6 |
| L9 | 418.98 | 395.23 | 389.81 | 350.68 | 383.05 | 382.18 | 371.21 | 361.22 | **359.47** | 396.17 | 440.82 | 432.3 | 426.23 | 418.65 | 381.1 | 381.38 | 381.85 | 367.78 | 387.66 | 385.9 | 389.47 | 385.02 |
| L10 | 469.5 | 393.79 | 385.79 | 351.87 | 394.39 | 375.03 | 367.68 | **361.98** | 364.83 | 389.09 | 427.62 | 438.39 | 433.79 | 436.38 | 369.76 | 369.61 | 371.04 | 370.45 | 372.08 | 373.18 | 388.66 | 376.47 |
| Average | 431.25 | 395.43 | 389.67 | 352.19 | 390.15 | 388.71 | **365.07** | 367.26 | 367.93 | 392.56 | 434.90 | 424.59 | 432.03 | 418.08 | 375.93 | 373.39 | 376.39 | 371.29 | 378.26 | 376.18 | 381.63 | 377.46 |

As for the comparison of neighborhood generation schemes, we can say that, the best scheme is insertion followed by mixed. Contrary to the conclusion of Tonyalı and Alkaya [40], 2-opt turned out to perform quite poorly so that for medium problems some of the algorithms perform even worse than CH.

Table 4. Per cent deviation of average scores from the optimum

| Small Problems | | | | Medium Problems | | |
|---|---|---|---|---|---|---|
| **Algorithm** | **Neighborhood** | **Dev. (%)** | | **Algorithm** | **Neighborhood** | **Dev. (%)** |
| MBO-DC-S | Insertion | 0.38 | | MBO-DC-S | Insertion | 1.88 |
| MBO-DC-M | Insertion | 0.46 | | SA-BSF | Insertion | 2.63 |
| MBO-DC-S | Mixed | 0.66 | | SA-BSF | Mixed | 3.1 |
| MBO-DC-M | Mixed | 0.69 | | MBO-IC-S | Insertion | 3.14 |
| SA-BSF | Mixed | 1.62 | | SA | Insertion | 3.17 |
| SA-BSF | Insertion | 1.65 | | MBO-DC-M | Mixed | 3.25 |
| MBO-IC-S | Insertion | 1.87 | | MBO-DC-S | Mixed | 3.26 |
| MBO-IC-S | Mixed | 1.89 | | SA | Mixed | 3.45 |
| SA | Insertion | 1.89 | | MBO-DC-M | Insertion | 3.76 |
| MBO-DC-S | 2-opt | 2.05 | | MBO-IC-S | Mixed | 4.5 |
| SA | Mixed | 2.05 | | MBO-IC-M | Insertion | 5.28 |
| MBO-IC-M | Mixed | 2.08 | | MBO-IC-M | Mixed | 5.47 |
| MBO-IC-M | Insertion | 2.22 | | SA | 2-opt | 10.44 |
| MBO-DC-M | 2-opt | 2.98 | | SA-BSF | 2-opt | 10.47 |
| MBO-IC-S | 2-opt | 5.72 | | MBO-DC-S | 2-opt | 11.35 |
| SA-BSF | 2-opt | 5.81 | | CHOrOpt | | 11.74 |
| MBO-IC-M | 2-opt | 6.62 | | MBO-DC-M | 2-opt | 12.24 |
| SA | 2-opt | 6.74 | | CH | | 13.31 |
| CHOrOpt | | 11.84 | | MBO-IC-S | 2-opt | 13.75 |
| CH | | 13.3 | | MBO-IC-M | 2-opt | 16.4 |
| NN | | 14.29 | | NN | | 22.74 |

The average run times of the algorithms (per run) are given in Table 5 where the runs are made on a Macbook Pro 2.3 GHz Dual-Core Intel Core i5, 8 GB RAM computer. The run times only slightly differed between algorithm variants and neighborhood functions and thus we preferred to display only the average figures here. We see that for the problem sizes undertaken in this study, it was possible to find the optimal solution by Gurobi in a reasonable time. We can also see that MBO consumes more time than SA which is inherent to its parallel work on multiple solutions in an iterative manner.

Table 5. Run times per problem per run in seconds.

| | Small | Medium | Large |
|---|---|---|---|
| NN | 0.00035 | 0.00170 | 0.00597 |
| CH | 0.04476 | 0.60539 | 8.80590 |
| CHOrOpt | 0.05796 | 0.93217 | 10.99521 |
| SA | 0.18430 | 7.97085 | 97.24959 |
| MBO | 1.16255 | 45.01208 | 1049.53680 |
| solver | 2.03487 | 8.32456 | 41.78564 |

We wanted to analyze what was more useful in MBO variants: S versus M or DC versus IC? For this, we wanted to refer Tables 1b and 1c since the most amount of statistics were collected for the small problems. From the average of averages, we obtain the results of M = 83.15, S = 82.82, DC = 82.10, IC = 83.87, from which we can conclude that switching from IC (immediate consideration of borrowed neighbors) to DC (delayed consideration of borrowed neighbors) is more effective than switching from M (multi step sharing of borrowed neighbors) to S (single step sharing of borrowed neighbors).

We would also like to say some words on the convergence and stagnancy behaviors of the SA and MBO algorithms. About the SA algorithm, for small problems and medium problems the BSF version helped a bit. On the contrary, for large problems, it seemed to be harmful. This can be interpreted as that for large problems the iteration limit was not sufficient and cutting the iterations at 90 per cent and going back to the best solution found so far acted as loss of time. As for the MBO, almost all solutions turned out to be the same for classical MBO (MBO-IC-M) on small problems. In 7500 neighbors they converged. When we switched to DC or S, convergence delayed almost to iteration limit K. For medium or large problems, we did not observe a clear convergence behavior neither for SA nor for MBO which is another indication of the insufficiency of the iteration limit for larger problems.

To summarize, we can claim that the suggested modifications on MBO (MBO-DC-S) has improved the performance of the original MBO (MBO-IC-M) by 1.84 and 3.40 percent on small and medium problems, averaging to 2.62 per cent. Also, we can say that MBO outperformed SA by 1.24 per cent on small problems and by 0.85 per cent on medium problems, averaging to 1.05 per cent.

## 6. Summary and conclusions

In this study, we undertook the MBO and SA algorithms and suggested some modifications on them related to their convergence behaviors. Then, we tested these modifications on a new application of the ATSP that is defined in this study. The extensive numerical experimentation reveals that, especially the modifications suggested for the MBO algorithm are quite successful and prevent premature convergence.

Two immediate continuations of this study can be explored in the future. First one is the case of multiple cooks and/or stoves in our application. Second one is the possibility of queueing up the raw meals in front of the stove while the cook prepares the next orders. This problem is much more complicated and will be sequence dependent.

**References**

[1] Jünger M, Reinelt G, and Rinaldi G (1995). The traveling salesman problem. *Handbooks in operations research and management science*, *7*, 225-330.

[2] Bläser M, Manthey B and Sgall J (2006). An improved approximation algorithm for the asymmetric TSP with strengthened triangle inequality. *Journal of Discrete Algorithms*, *4*(4), 623-632.

[3] Mömke T (2015). An improved approximation algorithm for the traveling salesman problem with relaxed triangle inequality. *Information Processing Letters*, *115*(11), 866-871.

[4] Roberti R and Toth P (2012). Models and algorithms for the asymmetric traveling salesman problem: an experimental comparison. *EURO J Transp Logist* **1**:113–133.

[5] Glover F, Gutin G, Yeo A and Zverovich A (2001). Construction heuristics for the asymmetric TSP. *European Journal of Operational Research*, *129*(3), 555-568.

[6] Duman E and Or I (2004). Precedence constrained TSP arising in printed circuit board assembly. *International Journal of Production Research*, *42*(1), 67-78.

[7] Subramanian A and Battarra M (2013). An iterated local search algorithm for the travelling salesman problem with pickups and deliveries. *Journal of the Operational Research Society*, *64*(3), 402-409.

[8] Uwaisy MA, Baizal ZKA and Reditya MY (2019). Recommendation of scheduling tourism routes using tabu search method (case study bandung). *Procedia Computer Science*, *157*, 150-159.

[9] Agrawal A, Ghune N, Prakash S and Ramteke M (2021). Evolutionary algorithm hybridized with local search and intelligent seeding for solving multi-objective Euclidian TSP. *Expert Systems with Applications*, *181*, 115192.

[10] Shi XH, Liang YC, Lee HP, Lu C and Wang QX (2007). Particle swarm optimization-based algorithms for TSP and generalized TSP. *Information processing letters*, *103*(5), 169-176.

[11] Skinderowicz R (2022). Improving Ant Colony Optimization efficiency for solving large TSP instances. *Applied Soft Computing*, *120*, 108653.

[12] Geng X, Chen Z, Yang W, Shi D and Zhao K (2011). Solving the traveling salesman problem based on an adaptive simulated annealing algorithm with greedy search. *Applied Soft Computing*, *11*(4), 3680-3689.

[13] Rao TS (2017). A Comparative Evaluation of GA and SA TSP in a Supply Chain Network. *Materials Today: Proceedings*, *4*(2), 2263-2268.

[14] Rao TS (2021). A simulated annealing approach to solve a multi traveling salesman problem in a FMCG company. *Materials Today: Proceedings*, *46*, 4971-4974.

[15] da Silva R, Venites Filho E and Alves A (2021). A thorough study of the performance of simulated annealing in the traveling salesman problem under correlated and long tailed spatial scenarios. *Physica A: Statistical Mechanics and its Applications*, *577*, 126067.

[16] Duman E, Uysal M and Alkaya AF (2012). Migrating birds optimization: a new metaheuristic approach and its performance on quadratic assignment problem. *Information Sciences*, *217*, 65-77.

[17] Benkalai I, Rebaine D, Gagné C and Baptiste P (2016). The migrating birds optimization metaheuristic for the permutation flow shop with sequence dependent setup times. *IFAC-PapersOnLine*, *49*(12), 408-413.

[18] Benkalai I, Rebaine D, Gagné C and Baptiste P (2017). Improving the migrating birds optimization metaheuristic for the permutation flow shop with sequence-dependent set-up times. *International Journal of Production Research*, *55*(20), 6145-6157.

[19] Meng T, Pan QK, Li JQ and Sang HY (2018). An improved migrating birds optimization for an integrated lot-streaming flow shop scheduling problem. *Swarm and Evolutionary Computation*, *38*, 64-78.

[20] Sioud A and Gagné C (2018). Enhanced migrating birds optimization algorithm for the permutation flow shop problem with sequence dependent setup times. *European Journal of Operational Research*, *264*(1), 66-73.

[21] Han Y, Li JQ, Gong D and Sang H (2018). Multi-objective migrating birds optimization algorithm for stochastic lot-streaming flow shop scheduling with blocking. *IEEE Access*, *7*, 5946-5962.

[22] Ping W, Sang HY, Tao QY and Qun S (2020, May). Improved Migratory Birds Optimisation Algorithm to Solve Low-Carbon Hybrid Lot-Streaming Flowshop Scheduling Problem. In *2020 World Conference on Computing and Communication Technologies (WCCCT)* (pp. 33-36). IEEE.

[23] Wang P, Sang H, Tao Q, Guo H, Li J, Gao K and Han Y (2020). Improved migrating birds optimization algorithm to solve hybrid flowshop scheduling problem with lot-streaming. *IEEE Access*, *8*, 89782-89792.

[24] Deng G, Xu M, Zhang S, Jiang T and Su Q (2022). Migrating birds optimization with a diversified mechanism for blocking flow shops to minimize idle and blocking time. *Applied Soft Computing*, *114*, 107834.

[25] Zhang S, Gu X and Zhou F (2020). An improved discrete migrating birds optimization algorithm for the no-wait flow shop scheduling problem. *IEEE access*, *8*, 99380-99392.

[26] Gao L and Pan QK (2016). A shuffled multi-swarm micro-migrating birds optimizer for a multi-resource-constrained flexible job shop scheduling problem. *Information Sciences*, *372*, 655-676.

[27] Zhang M, Tan Y, Zhu J, Chen Y and Chen Z (2020). A competitive and cooperative Migrating Birds Optimization algorithm for vary-sized batch splitting scheduling problem of flexible Job-Shop with setup time. *Simulation Modelling Practice and Theory*, *100*, 102065.

[28] Xiao Q, Guo X and Li D (2021). Partial disassembly line balancing under uncertainty: robust optimisation models and an improved migrating birds optimisation algorithm. *International Journal of Production Research*, *59*(10), 2977-2995.

[29] Zikai ZHANG, Qiuhua T, Zixiang L and Dayong H (2021). An efficient migrating birds optimization algorithm with idle time reduction for Type-I multi-manned assembly line balancing problem. *Journal of Systems Engineering and Electronics*, *32*(2), 286-296.

[30] Zhang Z, Tang Q, Han D and Li Z (2022). Multi-manned assembly line balancing with sequence-dependent set-up times using an enhanced migrating birds optimization algorithm. *Engineering Optimization*, 1-20.

[31] Ulker E and Tongur V (2017). Migrating birds optimization (MBO) algorithm to solve knapsack problem. *Procedia computer science*, *111*, 71-76.

[32] Aboudi NE and Benhlima L (2018, October). Towards parallel migrating birds framework for feature subset problem. In *Proceedings of the First International Conference on Data Science, E-learning and Information Systems* (pp. 1-5).

[33] Makas H and Yumusak N (2013, November). New cooperative and modified variants of the migrating birds optimization algorithm. In *2013 International Conference on Electronics, Computer and Computation (ICECCO)* (pp. 176-179). IEEE.

[34] Oz D (2017). An improvement on the Migrating Birds Optimization with a problem-specific neighboring function for the multi-objective task allocation problem. *Expert Systems with Applications*, *67*, 304-311.

[35] Niroomand S, Hadi-Vencheh A, Şahin R and Vizvári B (2015). Modified migrating birds optimization algorithm for closed loop layout with exact distances in flexible manufacturing systems. *Expert Systems with Applications*, *42*(19), 6586-6597.

[36] Cao J, Guan Z, Yue L, Ullah S and Sherwani RAK (2020). A bottleneck degree-based migrating birds optimization algorithm for the PCB production scheduling. *IEEE Access*, *8*, 209579-209593.

[37] Taşpınar N and Şimşir Ş (2020). An efficient SLM technique based on migrating birds optimization algorithm with cyclic bit flipping mechanism for PAPR reduction in UFMC waveform. *Physical Communication*, *43*, 101225.

[38] Tongur V, Ertunc E and Uyan M (2020). Use of the Migrating Birds Optimization (MBO) Algorithm in solving land distribution problem. *Land Use Policy*, *94*, 104550.

[39] Tongur V and Ülker E (2016). The analysis of migrating birds optimization algorithm with neighborhood operator on traveling salesman problem. In *Intelligent and Evolutionary Systems* (pp. 227-237). Springer, Cham.

[40] Tonyali S and Alkaya AF (2015). Application of Recently Proposed Metaheuristics to the Sequence Dependent TSP. In *Advanced Computational Methods for Knowledge Engineering* (pp. 83-94). Springer, Cham.

[41] A.F. Alkaya, E. Duman, "An Application of the Sequence Dependent Traveling Salesman Problem in Printed Circuit Board Assembly," *Transactions on Components, Packaging and Manufacturing Technology*, Vol.3, Issue 6, pp. 1063-1076, 2013.

[42] Babin G, Deneault S and Laporte G (2007). Improvements to the Or-opt heuristic for the symmetric travelling salesman problem. *Journal of the Operational Research Society*, *58*(3), 402-407.